\documentclass[12pt,thmsa]{article}
\usepackage{amsfonts}
\usepackage{amssymb}
\newtheorem{teo}{Theorem}[section]

\newtheorem{obs2}[teo]{Remark}

\newtheorem{tea}{Theorem}[subsection]

\newtheorem{no2}[teo]{Note}

\newtheorem{no3}[tea]{Note}

\newcommand{\Gal}{{\rm Gal}}
\newcommand{\Frob}{{\rm Frob }}

\newcommand{\mod}{{\rm mod}}

\newcommand{\Q}{\mathbb{Q}}

\newcommand{\F}{{\mathbb{F}}}

\begin{document}
\title{{\bf Elliptic $\mod \; \ell$ Galois representations which are not minimally elliptic
 }}

\author{Luis Dieulefait
\\
Dept. d'{\`A}lgebra i Geometria, Universitat de Barcelona;\\
Gran Via de les Corts Catalanes 585;
08007 - Barcelona; Spain.\\
e-mail: ldieulefait@ub.edu\\
 }
\date{\empty}

\maketitle

\vskip -20mm

\begin{abstract}
In a recent preprint (see [C]), F. Calegari has
shown that for $\ell = 2, 3, 5$ and $7$ there exist
 $2$-dimensional irreducible representations   $\rho$
 of $\Gal(\bar{\Q}/\Q)$ with values in
$\F_\ell$ coming from the $\ell$-torsion points of an
 elliptic curve defined over $\Q$, but not minimally,
  i.e., so that any elliptic curve giving rise to
  $\rho$ has  prime-to-$\ell$ conductor greater
   than the (prime-to-$\ell$) conductor of $\rho$.
In this brief note, we will show that the same is true
for any prime $\ell >7$.
\end{abstract}

\section{The result and its proof}

In this article, we are going to prove the following result:

 \begin{teo} For any prime $\ell >7$ the
 Galois representation $\rho$ obtained from the $\ell$ torsion points of the elliptic curve
 $$E^\ell: \qquad Y^2 = X (X- 3^\ell ) (X - 3^\ell - 1) $$
 is irreducible and unramified at $3$, but $E^\ell$ and any other elliptic curve
 giving rise to the same $\mod \; \ell$ Galois representation have bad reduction at $3$. Thus, these representations arise from 
  elliptic curves but not minimally. \\
   On the other hand, if we consider a modular abelian variety
   $A_f$ with good reduction at $3$ also giving rise to $\rho$
   (it follows from the modularity of elliptic curves and lowering the level that such a variety
    always exists) then as $\ell$ varies the dimension of
    $A_f$ tends to infinity with $\ell$.
    \end{teo}
    
  It was shown in [C] that also for primes $ \ell < 11 $  there exist irreducible residual representations that arise from elliptic curves but not minimally, thus the property is true for every prime.\\

 We will show that for every $\ell > 7$ the curve $E^\ell$
is semistable outside $2$, has bad reduction at $3$, the associated
 $\mod \; \ell$ Galois representation $\rho$ is
  irreducible, unramified at
$3$, and there is no elliptic curve with good
 reduction at $3$ whose associated $\mod \; \ell$ representation is isomorphic to $\rho$.\\

In page $9$ of [C], the example for $\ell = 7$ is constructed from the elliptic
curve $E$: 
$$ y^2 + y x + y = x^3 - 89 x + 316 $$
which has  semistable reduction at $2$ and its discriminant $\Delta$
 has $2$-adic valuation equal to $7$. This implies that the $\mod \; 7$ Galois
 representation $\rho$ attached to $E$ is unramified at $2$ (the prime-to-$7$ part of its conductor is $55$) and by Tate's theory satisfies:
 $a_2 \equiv \pm 3 \pmod{7}$, where $a_2$ is the trace of $\rho(\Frob \; 2)$.
 The representation can not correspond to an elliptic curve with good reduction
 at $2$ because for such an elliptic curve $E'$ we have $c_2 = 0, \pm 1, \pm 2$,
 where $c_2$ denotes the trace of the image of $\Frob \; 2$ for the compatible
 family of Galois representations attached to $E$, and therefore we would get
 $\pm 3 \equiv a_2 \equiv c_2 \pmod{7}$, a contradiction.\\

The same argument proves the result for higher primes: take $\ell
>7$ and consider the elliptic curve $E^\ell$. From the definition of $E^\ell$ we see that it has bad
 reduction at $2$ and $3$ and good reduction at $5$ and $\ell$.\\
 The same argument used for the case of the
Frey-Hellegouarch curves related to Fermat's Last Theorem (cf.
[H], pags. 368-369) shows that this curve is semistable outside
$2$ (i.e., it has semistable reduction at every odd prime of bad
reduction) and  that the
  corresponding $\mod \; \ell$ representation $\rho$ is
  unramified at $3$ (because the $3$-adic valuation of
  the minimal discriminant is multiple of $\ell$). From this and the fact that $E^\ell$ has bad semistable reduction
   at $3$ it follows that: $a_3 \equiv \pm 4 \pmod{\ell}$. \\

It is easy to check that $\rho$ is irreducible: in fact,
 this follows from the fact that it is semistable outside $2$ and
 has good reduction at $5$ (and comes from an elliptic curve).
 We indicate a short proof for the reader convenience: assume
 that $\rho$ is reducible, then (after semisimplifying, if
 necessary) we get: $\rho \cong \epsilon \oplus \epsilon^{-1} \chi \quad (*)$, where $\chi$
 denotes the $\mod \; \ell$ cyclotomic character and $\epsilon$ is a
 character unramified outside $2$ (here we use
  semistability outside $2$). Evaluating at $\Frob \; 5$ and taking traces we get:
   $a_5 \equiv r + 5 r^{-1} \pmod{\ell} \quad (**)$, where $ r = \epsilon
   (5)$. Since the $2$-part of the conductor of any elliptic curve
   is known to be at most $256$ it follows from (*) that the conductor of $\epsilon$ is at
   most $16$. Thus, since the image of $\epsilon$ is cyclic
   (it is contained in the multiplicative group of a finite field)
   we conclude that the order of $\epsilon$ is at most $4$, and in
   particular that $r^4 \equiv 1 \pmod{\ell}$.\\
On the other hand, we know that $a'_5 = 0, \pm 1, \pm 2, \pm 3,
\pm 4$, where $a'_5$ denotes the trace of the image of $\Frob \; 5$ for the $\ell$-adic representation attached to $E^\ell$,
 and thus $(a'_5 \; \mod \ell) = a_5$. With these restrictions
on $r$ and $a_5$, we can solve (**): squaring both sides
and using $r^2 = \pm 1$ and the above list of values for $a'_5$, we
check that the only possibility for (**) to hold is, if we
restrict to $\ell \geq 11$, $\ell = 17$ with $a'_5 = \pm 1$ (@).\\
This proves irreducibility for every $\ell \geq 11$, except for
$\ell = 17$. To rescue this last prime, observe that if we take
the curve $E^{17}$ we can count its number of points modulo $5$:
it has $8$ points. This gives $a'_5 = -2$ for $\ell = 17$. Hence,
since $-2 \neq \pm 1$, the case (@) never happens, and we also get
irreducibility for $\ell = 17$.\\

Since $a_3 \equiv \pm 4 \pmod{\ell}$ and $\ell \geq 11$, it is
clear that
 this representation can not correspond to an elliptic curve unramified
 at $3$, because for such an elliptic curve the corresponding trace
  $c_3$ at $\Frob \; 3$ ( in characteristic $0$)
satisfies $c_3 = 0, \pm 1, \pm 2, \pm 3$, thus $a_3 \equiv c_3 \pmod{\ell}$
gives a contradiction.\\

Since all elliptic curves over $\Q$ are modular, by level-lowering
we know that there exists a weight $2$ newform $f$ of level prime to $3$
(and equal to the prime-to-$\ell$ part of the conductor of $\rho$) such
that $\ell$ splits totally in the field
 $\Q_f$ generated by the eigenvalues of $f$ and for a prime
  $\lambda \mid \ell$ in $\Q_f$ the $\mod \; \lambda$ representation
   $\bar{\rho}_{f,\lambda}$ attached to $f$ is isomorphic to $\rho$.
    Of course, due to the result proved above, it must hold
 $\Q_f \neq \Q$, so that the abelian variety $A_f$ associated to $f$
  is not an elliptic curve. \\
 Moreover, it is not hard to see that given any dimension $d$, for
 almost every prime $\ell$ any abelian variety $A_f$ realizing $\rho$
 with minimal ramification as above (i.e., with the level of $f$ equal
  to the prime-to-$\ell$ part of the conductor of $\rho$ and the
  residual representation attached to $f$ isomorphic to $\rho$)
   must be of dimension grater than $d$. This follows from the
   fact that if the dimension is bounded by $d$, the degree of
    the field generated by $c_3$, the trace at $\Frob \; 3$ of the Galois
     representations attached to $f$, is also bounded by $d$, and from
     this it follows (using the bound for the complex absolute values
     of $c_3$ and its Galois conjugates) that there are only finitely
     many possible values for $c_3$. Since (again) $c_3 \neq \pm 4$,
 the congruence $a_3 \equiv c_3$ gives
 $$ c_3 \equiv \pm 4 \pmod{\ell}$$
 which can only be satisfied by finitely many primes $\ell$ (for a fixed $d$), and this
  is what we wanted to prove.\\

\section{Bibliography}

[C] Calegari, F., {\it Mod $p$ representations on Elliptic Curves}, preprint, available at\\
http://front.math.ucdavis.edu/math.NT/0406244
\\

[H] Hellegouarch, Y., {\it Invitation to the Mathematics of
Fermat-Wiles}, Academic Press, 2002

\end{document}